\documentclass[11pt]{article}
\usepackage{amssymb, amsfonts, amsmath }
\usepackage{latexsym}

\providecommand{\MR}{\relax\ifhmode\unskip\space\fi MR }

\def\Fix{{\it Fix \,}}
\def\Per{{\it Per \,}}

\newtheorem{thm}{Theorem}[section]
\newtheorem{prop}[thm]{Proposition}
\newtheorem{lem}[thm]{Lemma}

\newcommand{\demo}{ \noindent {\it   Proof. }}

\title{Free-by-cyclic groups have solvable conjugacy problem}

\author{O.\ Bogopolski \\ \small{Institute of Mathematics of} \\ \small{Siberian Branch of Russian Academy
of Sciences,} \\ {\small Novosibirsk, Russia} \\ \small{e-mail:
groups@math.nsc.ru}\\ \\ \\ A.\ Martino \\ \small{Centre de
Recerca Matematica} \\ \small{Bellaterra, Spain} \\
\small{e-mail: AMartino@crm.es} \\ \\ \\ O.\ Maslakova
\\ \small{Institute of Mathematics of} \\ \small{Siberian Branch of Russian Academy
of Sciences,} \\ {\small Novosibirsk, Russia} \\ \small{e-mail: tessae@ngs.ru}\\ \\ \\ E.\ Ventura \\
\small{Dept.\ Mat.\ Apl.\ III, Univ.\ Pol.\ Catalunya, Barcelona,
Spain}
\\ \small{and Dept.\ of Math.,\, University of Nebraska-Lincoln}
\\ \small{e-mail: enric.ventura@upc.es} }

\begin{document}

\maketitle

\begin{abstract}
We show that the conjugacy problem is solvable in [finitely generated free]-by-cyclic groups, by using a result of \,
O. Maslakova that one can algorithmically find generating sets for the fixed subgroups of free group automorphisms, and
one of P. Brinkmann that one can determine whether two cyclic words in a free group are mapped to each other by some
power of a given automorphism. The algorithm effectively computes a conjugating element, if it exists. We also solve
the power conjugacy problem and give an algorithm to recognize if two given elements of a finitely generated free group
are Reidemeister equivalent with respect to a given automorphism.
\end{abstract}

\section{Introduction}

A \emph{free-by-cyclic} group is a group $G$ having a free normal
subgroup $F$ with cyclic quotient $C=G/F$. If $F$ can be chosen to
be finitely generated then $G$ is a \emph{[f.g.\,free]-by-cyclic}
group. (Note that the parenthesis are relevant here since surface
groups are both free-by-cyclic and finitely generated, but most of
them are not [f.g.\,free]-by-cyclic.)

We shall be concerned with [f.g.\,free]-by-cyclic groups. Let
$x_1, \ldots , x_n$ be a basis for $F$ and let $t$ be a pre-image
in $G$ of a generator of $C$. Right conjugation by $t$ in $G$
induces an automorphism of $F$, which we denote $\phi$. Note that
we shall write $\phi$ on the right so that the image of $w$ under
$\phi$ will be denoted $w \phi$.

The following is then a presentation of $G$
 $$
G=\langle x_1, \ldots ,x_n ,t \; :\;  t^{-1} x_i t = x_i \phi ,\; t^m=h \rangle,
 $$
where $m$ is the cardinal of $C$ and $h$ is some element in $F$,
understanding that the relation $t^m=h$ is not present when
$m=\infty$. Our main result is the following

\begin{thm}\label{main}
The conjugacy problem in [f.g.\,free]-by-cyclic groups is solvable.
\end{thm}

For some special cases, this result is already known. The
automorphism $\phi$ is said to have no periodic conjugacy classes
if one cannot find an integer $k$ and elements $g,h \in F$ such
that $g \phi^k = h^{-1} gh$. If $C$ is infinite, this is
equivalent to saying that $G$ has no $\mathbb{Z} \oplus
\mathbb{Z}$ subgroups and hence, by \cite{brinkmann}, \cite{BF}
and \cite{BF2}, that $G$ is hyperbolic. Also, if $C$ is finite,
the group $G$ is virtually free and hence hyperbolic. In all these
cases, then, it is well known that $G$ has solvable conjugacy
problem. But, clearly, not all [f.g.\,free]-by-cyclic groups are
hyperbolic. It has even been announced, in \cite{BR}, that they
fail to be automatic, in general.

Some other partial results are already known in this direction.
For example, in the preprint~\cite{BBV} the authors consider the
case where some power of $\phi$ is an inner automorphism, and give
an algorithm to decide if two given elements in $G$ of the form
$tu,\, tv$ for $u,v\in F$ are conjugated by some element in $F$.

Our proof of Theorem~\ref{main} will work, in general, for arbitrary [f.g.\,free]-by-cyclic groups, including the
previously known particular cases. The algorithm provided also computes a conjugating element, when it exists. Our
proof relies on the following recent theorems.

\begin{thm}[Maslakova, \cite{Mas}]
\label{mas} There exists an algorithm to compute a finite generating set for the fixed point subgroup of an arbitrary
automorphism of a free group of finite rank.
\end{thm}

\begin{thm}[Brinkmann, \cite{brinkmann2}]
\label{map} Given a finitely generated free group $F$, two elements $u,v$ of $F$ and an automorphism $\phi$ of $F$ it
is decidable whether there exists an integer $k$ such that $u \phi^k$ is conjugate to $v$. Moreover, if such a $k$
exists, it can be computed.
\end{thm}

Our solution to the conjugacy problem for $G$ proceeds by showing that, in the light of Theorem~\ref{map}, it can be
reduced to the \emph{twisted conjugacy problem} for $F$. Then, we solve this classical problem providing an algorithm
to recognize Reidemeister classes with respect to automorphisms of finitely generated free groups.

Let $F$ be a free group and $\phi$ an automorphism of $F$. Two elements $u,v \in F$ are said to be \emph{$\phi$-twisted
conjugate}, denoted $u\sim_{\phi} v$, if there exists $g\in F$ such that $(g\phi)^{-1} u g=v$. The equivalence relation
$\sim_{\phi}$ was first introduced by Reidemeister in~\cite{Reid}, and has an important role in Nielsen fixed point
theory. A couple of interesting references are~\cite{LeLu}, where it is proven that the number of $\phi$-twisted
conjugacy classes is always infinite, and~\cite{DV}, where Problem~3(i) in the Open Problem section asks for an
algorithm recognizing $\phi$-twisted conjugacy classes.

It is said that the \emph{$\phi$-twisted conjugacy problem} is solvable in $F$ if, for any elements, $u,v\in F$, we can
algorithmically decide if $u\sim_{\phi} v$ (for example, the $\operatorname{id}$-twisted conjugacy problem is the
standard conjugacy problem in $F$, which is clearly solvable). And it is said that the \emph{twisted conjugacy problem}
is solvable in $F$ if the $\phi$-twisted conjugacy problem is solvable for any $\phi \in Aut(F)$. This twisted
conjugacy problem is also part of a more general problem posted by G. Makanin in Question~10.26(a) of~\cite{Mak}.

We first prove the following.

\begin{prop}\label{red}
Let $F$ be a finitely generated free group. If the twisted conjugacy problem is solvable in $F$ then the standard
conjugacy problem is solvable in $G$.
\end{prop}

Then, we give a solution for the twisted conjugacy problem in a finitely generated free group $F$, thus extending the
result in~\cite{BBV}, where it is done for virtually inner automorphisms, and answering question 3(i) of~\cite{DV} in
the affirmative.

\begin{thm}\label{tcp}
Let $F$ be a finitely generated free group. The twisted conjugacy problem is solvable in $F$.
\end{thm}

Now, Theorem~\ref{main} follows immediately from Proposition~\ref{red} and Theorem~\ref{tcp}.

Finally, in the last section we develop few technical lemmas that will allow us to extend Theorem~\ref{main}
(essentially with the same proof) to the following result.

\begin{thm}\label{pcp}
The power conjugacy problem in [f.g. free]-by-cyclic groups is solvable.
\end{thm}

In an arbitrary group, two elements $u$ and $v$ are said to be \emph{power conjugated} when there exist integers $p,q$
such that $u^p$ and $v^q$ are non-trivial and conjugated to each other in the group. The \emph{power conjugacy problem}
in a group consists on deciding whether two given elements are power conjugated (and find such exponents and
conjugating element if they exist).

\medskip

Before going into the details of the algorithm, we make some
remarks.

M. Lustig has a recent series of two preprints, \cite{Lu},
indicating a solution to the conjugacy problem in $Aut(F)$ and
$Out(F)$. As a consequence, he obtains also an algorithm for
computing the fixed subgroup of any automorphism of $F$, thus
providing an alternative proof for Theorem~\ref{mas}. It also
seems that these preprints implicitly contain a solution for the
twisted conjugacy problem in $F$.

It was pointed out to us by Ilya Kapovich that many one-relator groups are [f.g.\,free]-by-cyclic (one can prove this
by imposing few assumptions on the relator). It seems that there is evidence to think that these assumptions are quite
weak, meaning that most relators satisfy them and, hence, most of 1-relator groups fall in the family of groups
considered in this paper. However, this has not been expressed yet in a precise form at this time.

It is worth mentioning that in Chapter~3 of~\cite{Miller} there is an explicit construction of a [f.g.\,
free]-by-[f.g.\, free] group with unsolvable conjugacy problem. So, our result is no longer true if we replace
``cyclic" by an arbitrary ``finitely generated free group".

Finally, we state the following lemma for later use.

\begin{lem}\label{twclasses}
Let $\phi$ be an automorphism of a free group $F$. Then, any $\phi$-twisted conjugacy class in $F$ is a union
of $\phi$-orbits.
\end{lem}

\demo It is sufficient to prove that if two elements from $F$ lie in the same $\phi$-orbit, then they lie in the same
$\phi$-twisted conjugacy class. By induction, this reduces to prove that, for every $u\in F$, $u\sim_{\phi} u\phi$. And
this fact is obvious since $u=(u\phi)^{-1} (u\phi) u$. $\Box$

\section{The conjugacy problem}

First note that, using the relations $wt=t(w\phi)$ and $wt^{-1}=t^{-1}(w\phi^{-1})$ for $w\in F$, every element in $G$
can be algorithmically re-written as a word of the form $t^r u$, where $r$ is an integer and $u\in F$. In the case
where $C$ is finite and we also have the relation $t^m=h$, we may further assume that $0\leq r\leq m-1$. In either
case, we get a unique representation for elements of $G$, which is algorithmically computable from a given arbitrary
word on the generators.

If we conjugate $t^r u$ by an arbitrary element $t^k g$, we obtain
 $$
(t^k g)^{-1}(t^r u)(t^k g)=t^r (g\phi^r)^{-1} t^{-k} ut^k g = t^r (g\phi^r)^{-1} (u \phi^k )g.
 $$
Hence, two elements in $G$, say $t^r u$ and $t^s v$ (with $0\leq r,s\leq m-1$ in the case where $|C|=m<\infty$), are
conjugate in $G$ if and only if $r=s$ and $v\sim_{\phi^r} (u\phi^k )$ for some integer $k$. This is the key fact in the
following discussion.

\bigskip

\noindent \emph{Proof of Proposition~\ref{red}.} Suppose two elements in $G$ are given, say $t^r u$ and $t^s v$. We
have to decide if they are conjugate to each other in $G$, and find a conjugating element if it exists.

We first deal with the case where $r=0$. Note that $u$ is only conjugate in $G$ to other elements $v$ of the base group
$F$. Moreover, $u$ is conjugate to $v$ in $G$ if and only if, some power of the automorphism $\phi$, maps $u$ to a
conjugate of $v$. This is decidable by Theorem~\ref{map}, so we can decide if $u,v\in F$ are conjugate in $G$.

For the case $r\neq 0$ note that, by Lemma~\ref{twclasses}, $u \phi^k \sim_{\phi^r} u\phi^{k\pm r}$. Hence, $t^r u$ and
$t^s v$ are one conjugate to the other in $G$ if, and only if, $r=s$ and $v \sim_{\phi^r} (u\phi^k )$ for some integer
$0\leq k\leq |r|-1$. Thus, a solution for the twisted conjugacy problem in $F$ provides a solution for the standard
conjugacy problem in $G$. $\Box$

\bigskip

\noindent \emph{Proof of Theorem~\ref{tcp}.} Let $\phi$ be an automorphism of $F$, and suppose $u,v \in F$ are given.
We need to algorithmically decide whether $u\sim_{\phi} v$.

Choose a free basis for $F$ and, adding a new letter $z$, we get a free basis for $F\,'=F* \langle z \rangle$. Let
$\phi'\in Aut(F\,')$ be the extension of $\phi$ defined by $z\phi' =uzu^{-1}$. Let $\gamma_y$ denote the inner
automorphism of $F\,'$ given by right conjugation by $y\in F\,'$, $x \gamma_y =y^{-1}xy$.

We claim that $u\sim_{\phi} v$ if, and only if, $\Fix{(\phi'\gamma_v)}$ contains an element of the form $g^{-1}zg$ for
some $g\in F$ (and, in this case, $g$ itself is a valid $\phi$-twisted conjugating element).

In fact, suppose that $v=(g\phi)^{-1} ug$ for some $g\in F$. A simple calculation shows that $g^{-1}zg$ is then fixed
by $\phi' \gamma_v$. Conversely, if $g^{-1}zg$ is fixed by $\phi' \gamma_v$ for some $g\in F$, then
$gv^{-1}(g\phi)^{-1}u$ commutes with $z$. And this implies $gv^{-1}(g\phi)^{-1}u=1$, since this word contains no
occurrences of $z$. Hence, $v=(g\phi)^{-1}ug$ and $u\sim_{\phi} v$ (with $g$ being a $\phi$-twisted conjugating
element).

Since, by Theorem~\ref{mas}, we can algorithmically find a
generating set for $\Fix{(\phi' \gamma_v)}$, we can also decide if
this subgroup contains an element of the form $g^{-1}zg$ for some
$g\in F$. One can, for example, look at the corresponding (finite)
core-graph for $\Fix{(\phi' \gamma_v)}$ (algorithmically
computable from a set of generators) and see if there is some loop
labelled $z$ at some vertex connected to the base-point by a path
whose label does not use the letter $z$. If this is the case, the
label of such a path provides the $g$, i.e. the required
$\phi$-twisted conjugating element.

(It is not difficult to show that $\Fix{(\phi' \gamma_v)}$
contains an element of the form $g^{-1}zg$ if, and only if, it
contains some word involving the letter $z$; and, in this case,
the longest initial $F$-segment in such a word provides the
$\phi$-twisted conjugating element. With this observation, one can
slightly simplify the algorithm given, by just checking to see
whether any of the generators of $\Fix{(\phi' \gamma_v)}$ involve
$z$.) $\Box$

\section{The power conjugacy problem}

With the help of a few technical lemmas, the argument given to solve the conjugacy problem in $G$ also works, in much
the same way, to solve the power conjugacy problem.

Theorem~\ref{mas} can be extended to consider periodic subgroups.
Recall that, given an automorphism $\phi$ of $F$, the
\emph{periodic subgroup} of $\phi$ is the subgroup $$\Per{\phi}=\{
w \in F : w \phi^k = w \ \mbox{\rm for some } k>0
\}=\cup_{k=1}^{\infty} Fix\, \phi^k.$$

\begin{prop}
\label{per} There exists an algorithm to compute a finite generating set for the periodic subgroup of any given
automorphism, $\phi$, of a finitely generated free group $F$. More precisely, there exists a computable integer $p_0$
(independent of $\phi$) such that $\Per{\phi}=\Fix{\phi^{p_0}}$.
\end{prop}

\demo It is well known that, for a finitely generated free group $F$ of rank $n\geq 0$, the group $Aut(F)$ has bounded
torsion (Stallings first proved this in~\cite{St}). What we need here is a computable integer $p_0$ (only depending on
$n$) such that the order of any finite order element in $Aut(F)$ divides $p_0$. A possible direct proof follows
(see~\cite{LeNi} and~\cite{K} for better bounds that can possibly reduce the complexity of our algorithm).

Clearly, if $n=0$ or $n=1$, we can take $p_0 =2$. For $n\geq 2$, we can invoke Theorem~2.1 of~\cite{Cul}, which implies
that every finite order element of $Out(F)$ can be realised as a graph automorphism of a finite graph $Z$ with rank
$n$. Deleting the degree 1 and degree 2 vertices in $Z$, we can assume that $Z$ contains no such vertices. It is easy
to see then that $Z$ has at most $3n-3$ edges, which total to a maximum of $6n-6$ oriented edges. Hence, every finite
order element of $Out(F)$ has order dividing $p_0 =(6n-6)!$. The same is true for $Aut(F)$, since the natural map
$Aut(F) \to Out(F)$ has torsion free kernel.

The result we have to prove is clear when $F$ is cyclic. So, we can assume $n\geq 2$.

By Corollaries~3.6 and~3.7 of \cite{St}, there exists $s\geq 1$
such that $\Per{\phi}=\Fix{\phi^{s}}$ (we assume further that $s$
is minimal possible). In particular, $\Per{\phi}$ has rank $r\leq
n$ (see~\cite{BH}), and $\phi$ restricts to an automorphism
$\phi'\in Aut(\Per{\phi})$ of order $s$. So, $s$ either divides 2
(if $r=0,1$) or $(6r-6)!$ (otherwise). In any case, $s$ divides
$p_0 =(6n-6)!$. So, $\Per{\phi}=\Fix{\phi^{p_0}}$.

Finally, using Theorem~\ref{mas}, we are done. $\Box$

\bigskip

Let $\phi$ be an automorphism of a finitely generated free group $F$.

For any $p\geq 1$, and any $w\in F$, we define $w_{\phi,\, p}=(w
\phi^{p-1}) (w \phi^{p-2}) \cdots (w\phi)w$. This notation will be
useful because, for every integer $r$ and every $u\in F$, we have
$(t^r u)^p =t^{rp}u_{\phi^r,\,p}$ in $G$. Note that
$w_{\operatorname{id},\,p}=w^p$. Note also that, for every
automorphism $\psi$ commuting with $\phi$, we have $w_{\phi,\,
p}\psi =(w\psi)_{\phi,\, p}$. Also, $u\sim_{\phi} v$ implies
$u\psi \sim_{\phi} v\psi$.

\begin{lem}\label{tecnic}
Let $u,\,v\in F$. If $u\sim_{\phi} v$ then $u_{\phi,\, p}\sim_{\phi^p} v_{\phi,\, p}$ for every $p\geq 1$.
\end{lem}

\demo Assume the existence of an element $g\in F$ satisfying $(g\phi)^{-1} ug=v$. Then, applying $\phi^i$ on both
sides, we obtain $(g\phi^{i+1})^{-1}(u\phi^i)(g\phi^i)=v\phi^i$. Now, multiplying all these equations,
 $$
(g\phi^p)^{-1} u_{\phi,\, p} g=\Pi_{i=p-1}^{0} (g\phi^{i+1})^{-1}(u\phi^i)(g\phi^i) = \Pi_{i=p-1}^{0} v\phi^i =
v_{\phi,\, p}.
 $$
This proves that $u_{\phi,\, p}\sim_{\phi^p} v_{\phi,\, p}$ (with the same twisted conjugating element $g$). $\Box$

\bigskip

Adapting the proof of Theorem~\ref{tcp}, we can obtain the following technical result.

\begin{lem}\label{powertwist}
Given $u,v\in F$, one can algorithmically decide if $u_{\phi,\, p} \sim_{\phi^p} v_{\phi,\, p}$ for some $p\geq 1$.
\end{lem}

\demo As before, add a new generator $z$ to $F$ and consider the extension $\phi' \in Aut(F* \langle z \rangle )$ of
$\phi$ given by $z\phi' =uzu^{-1}$. Exactly the same arguments as above show now that, for $p\geq 1$ and $g\in F$,
$v_{\phi,\, p} =(g\phi^p)^{-1} u_{\phi,\, p} g$ if and only if $g^{-1}zg \in \Fix{(\phi' \gamma_{v})^p}$ (to do this
computation, note that $(\phi'\gamma_v)^p =\phi^{\prime p}\gamma_{v_{\phi,\,p}}$ and $z\phi^{\prime p}=u_{\phi,\,
p}zu_{\phi,\,p}^{-1}$). So, we are done by invoking Proposition~\ref{per} (and the computability of $p_0$ there ensures
that we can compute the value of $p$ here). $\Box$

\bigskip

Now, we can adapt the proof of Proposition~\ref{red} to solve the power conjugacy problem in $G$.

\bigskip

\noindent \emph{Proof of Theorem~\ref{pcp}.} Suppose we are given two elements, $t^r u$ and $t^s v$, from $G$,
$r,\,s\in \mathbb{Z}$, $u,v\in F$. We need to decide whether they are power conjugated in $G$. Note that, if these two
elements have infinite order, this is the same as deciding whether there exist non-zero exponents $p$ and $q$ such that
$(t^r u)^p$ and $(t^s v)^q$ are conjugated to each other in $G$.

As before, we deal first with the case $r=s=0$. Here, given
$u,v\in F$ we have to decide whether for some integers $p,q\neq
0$, $u^p$ is mapped to a conjugate of $v^q$ by some power of
$\phi$. In a free group we can algorithmically find {\em roots} of
elements. That is, there is a unique element $\hat{u}$ of $F$ such
that $u$ is a positive power of $\hat{u}$, and $\hat{u}$ is not
itself a proper power. Similarly, there exists a root $\hat{v}$,
for $v$. Thus there exist integers $k_1, k_2$ such that
$u^p=\hat{u}^{k_1}$ and $v^q=\hat{v}^{k_2}$.

Since roots are unique in free groups, $u^p$ is mapped to a
conjugate of $v^q$ by some power of $\phi$ if, and only if,
$\hat{u}$ is mapped to a conjugate of $\hat{v}^{\epsilon}$, for
$\epsilon={\pm 1}$ and $k_1=\epsilon k_2$. Thus, we are done by
invoking Theorem~\ref{map} applied to the roots and observing that
$k_1$ and $k_2$ are computable.

Now, if $|C|=m$ is finite, then any element of $G$ raised to the
power $m$ lies in $F$. Moreover, $t^r u$ has infinite order if and
only if, $(t^r u)^m \neq 1$. Thus, if $t^r u$ and $t^s v$ are to
be power conjugated, then either $(t^r u)^m=1=(t^s v)^m$ or $1
\neq (t^r u)^m, (t^s v)^m \in F$. In the former case we can solve
the power conjugacy problem by finitely many checks of the
standard conjugacy problem for $G$. In the latter case, we are
done by the argument in the preceding paragraph, applied to the
pair $(t^r u)^m, (t^s v)^m$.

So, we can restrict our attention to the case $m=\infty$. In particular, $G$ is torsion-free.

By applying the future algorithm twice (once for the pair of elements $t^r u$, $t^s v$, and again for $t^r u$, $(t^s
v)^{-1}$) we may restrict our attention to positive exponents, $p,q$. Note that, if for some integers $p,q\geq 1$,
$(t^r u)^p =t^{rp}u_{\phi^r,\, p}$ and $(t^s v)^q =t^{sq}v_{\phi^s,\, q}$ are conjugated to each other in $G$, then
$rp=sq$. In particular, $r=0$ if and only if $s=0$. And if both $r$ and $s$ are not zero then $t^r u$ and $t^s v$ are
power conjugated in $G$ if and only if $(t^r u)^s$ and $(t^s v)^r$ also are. Thus, our problem reduces to the case
$r=s$.

Since we have dealt with the case $r=s=0$ above, it remains to consider the situation where $r=s\neq 0$ (and hence,
$p=q$). That is, we are given elements of the form $t^r u$ and $t^r v$ with $r\neq 0$, and we have to decide if there
exists an integer $p\geq 1$ such that $(t^r u)^p$ and $(t^r v)^p$ are conjugate to each other in $G$, i.e. such that
$v_{\phi^r,\, p}\sim_{\phi^{rp}} (u_{\phi^r,\, p}\phi^k)$ for some integer $k$.

Here, we claim that $(u_{\phi^r,\, p}\phi^k) \sim_{\phi^{rp}} (u_{\phi^r,\, p}\phi^{k\pm r})$ for every $p\geq 1$. In
fact, by Lemma~\ref{twclasses}, we have $u\sim_{\phi^r} (u\phi^{\pm r})$ so, using Lemma~\ref{tecnic}, $u_{\phi^r,\, p}
\sim_{\phi^{rp}} (u\phi^{\pm r})_{\phi^r,\, p}=u_{\phi^r,\, p}\phi^{\pm r}$. Then, $(u_{\phi^r,\, p}\phi^k
)\sim_{\phi^{rp}} (u_{\phi^r,\, p}\phi^{k\pm r})$, for every $k,\, p\in \mathbb{Z}$, $p\geq 1$. Thus it only remains to
decide whether there exists an integer $p\geq 1$ such that $v_{\phi^r,\, p}\sim_{\phi^{rp}} (u_{\phi^r,\, p}\phi^k)$
for some integer $0\leq k\leq |r|-1$. But $(u_{\phi^r,\, p}\phi^k)=(u\phi^k)_{\phi^r,\, p}$ so, using
Lemma~\ref{powertwist} at most $r$ times, we are done. $\Box$

\section*{Acknowledgments}

We thank S. Hermiller, I. Kapovich, G. Levitt and M. Lustig for interesting comments on the subject. The first named
author is partially supported by the grant of the President of Russian Federation for young Doctors MD-326.2003.01, and
by the INTAS grant N 03-51-3663. The second named author gratefully acknowledges the postdoctoral grant SB2001-0128
funded by the Spanish government, and thanks the CRM for its hospitality during the academic course 2003-2004. The
third named author is partially supported by the Grant Council of the President of Russian Federation through grant
NS-2069.2003.1, and by Lavrent'ev's grant for young scientists of the Siberian Branch of the Russian Academy of
Sciences. The forth named author gratefully acknowledges partial support by DGI (Spain) through grant BFM2003-06613,
and by the Generalitat de Catalunya through grant ACI-013. He also thanks the Department of Mathematics of the
University of Nebraska-Lincoln for its hospitality during the second semester of the course 2003-2004, while this
research was conducted.

\bigskip

\bibliographystyle{amsplain}

\begin{thebibliography}{99}

\bibitem{BBV} V.~Bardakov, L.~Bokut and A.~Vesnin, \emph{Twisted conjugacy in free groups and Makanin's question},
preprint, http://arxiv.org/abs/math.GR/0401349.

\bibitem{BF} M.~Bestvina and M.~Feighn, \emph{A combination theorem for negatively curved groups}, J. Differential Geom.
\textbf{35} (1992), no.~1, 85--101.

\bibitem{BF2} M.~Bestvina and M.~Feighn, \emph{Addendum and correction to: ``{A} combination theorem for
negatively curved groups''}, J. Differential Geom. \textbf{43}
(1996), no.~4, 783--788.

\bibitem{BH} M. Bestvina, M. Handel, \emph{Train tracks and automorphisms of free groups}, Ann. of Math., {\bf 135} (1992),
1--51.

\bibitem{BR} M.R. Bridson and L.~Reeves, \emph{On the absence of automaticity of certain free-by-cyclic groups}, in
preparation.

\bibitem{brinkmann} P.~Brinkmann, \emph{Hyperbolic automorphisms of free groups}, Geom. Funct. Anal. \textbf{10} (2000),
no.~5, 1071--1089.

\bibitem{brinkmann2} P.~Brinkmann, \emph{Dynamics of free group automorphisms},
preprint. Available at http://uk.arxiv.org/abs/math.GR/0308199.

\bibitem{Cul} M. Culler, \emph{Finite groups of outer automorphisms of a free group}, Contributions to group theory,
197--207, Contemp. Math., 33, Amer. Math. Soc., Providence, RI, 1984.

\bibitem{DV} W. Dicks, E. Ventura, \emph{The group fixed by a family of injective endomorphism of a free group},
Contemp. Math., {\bf 195} (1996), 1-81.

\bibitem{K} D.~G. Khramtsov, \emph{Finite groups of automorphisms of free
groups}, Math. Notes, {\bf 38} (1985), 721-724; transl. from
Matem. Zametki, {\bf 38}, N 3 (1985), 386-392.

\bibitem{LeLu} G. Levitt, M. Lustig, \emph{Most automorphisms of a hyperbolic group have very simple dynamics}, Ann.
Scient \'{E}c. Norm. Sup., \textbf{33} (2000), 507-517.

\bibitem{LeNi} G. Levitt, J.L. Nicolas, \emph{On the maximum order of torsion elements in $GL(n,\mathbb{Z})$ and
$Aut(F_n)$}, Journal of Algebra, \textbf{208} (1998), 630-642.

\bibitem{Lu} M. Lustig, \emph{Structure and conjugacy for automorphisms of free groups I, II}, Max-Planck Institut f\"ur
Mathematik Preprint Series 2000, no. 130, and 2001, no. 4 (see
http://www.mpim-bonn.mpg.de).


\bibitem{Mas} O.~S. Maslakova, \emph{The fixed point group of a free group automorphism}, Algebra i Logika, \textbf{42} (2003), no.~4,
422--472. Translated (English) Algebra and Logic, \textbf{42}, no.
4, 2003, 237 - 265.


\bibitem{Miller} C.F. Miller, \emph{On group-theoretic decision problems and their classification}, Annals of Mathematics
Studies, No. 68. Princeton University Press, Princeton, N.J.; University of Tokyo Press, Tokyo, 1971. viii+106 pp.

\bibitem{Reid} K. Reidemeister, \emph{Automorphismen von homotopiekettenringen}, Math. Ann., \textbf{112} (1936),
586-593.

\bibitem{St} J.R. Stallings, \emph{Finiteness properties of matrix representations}, Annals of Mathematics,
\textbf{124} (1986), 337-346.

\bibitem{Mak} \emph{The Kourovka notebook: Unsolved problems in group theory} (Russian). Fifteenth augmented edition. Edited by V.
D. Mazurov and E. I. Khukhro. Russian Academy of Sciences,
Siberian Division, Institute of Mathematics, Novosibirsk, 2002.
172 pp.
\end{thebibliography}

\end{document}